\documentclass[a4paper,12pt]{article}

\usepackage[T2A]{fontenc}
\usepackage[utf8]{inputenc}   % older versions of ucs package
\usepackage[russian]{babel}

\usepackage{CJKutf8}

\usepackage{amsmath,amstext,amsfonts,amsthm,amssymb}
\usepackage{eucal}
\usepackage{graphicx}
\renewcommand{\subsubsection}[1]{\addtocounter{subsubsection}{1}
{\ \\[3pt]\bf \thesubsubsection. \  #1} }
\swapnumbers

{  \theoremstyle{definition}

}
\newcommand{\chr}{{\operatorname{ch}}}

\newcommand{\iso}{\overset{\sim}{\longrightarrow}}
\newcommand{\isom}{\overset{\sim}{=}}
\newcommand{\lra}{\longrightarrow}

\newcommand{\dpar}{\partial}

\newcommand{\bea}{\begin{eqnarray*}}
\newcommand{\eea}{\end{eqnarray*}}
\newcommand{\bean}{\begin{eqnarray}}
\newcommand{\eean}{\end{eqnarray}}

\newcommand{\fl}{\mathfrak l}

\newcommand{\CA}{\mathcal{A}}

\newcommand{\CC}{\mathcal{C}}

\newcommand{\CD}{\mathcal{D}}

\newcommand{\CK}{\mathcal{K}}

\newcommand{\BC}{\mathbb{C}}

\newcommand{\BP}{\mathbb{P}}
\newcommand{\BQ}{\mathbb{Q}}

\newcommand{\BZ}{\mathbb{Z}}

\newcommand{\nc}{\newcommand}

\nc{\Id}{\text{Id}}
\nc{\la}{\lambda}

\begin{document}

%\centerline{\b }

\bigskip\bigskip

%\centerline{\bf INFINITE SPHERES AND KOSZUL COMPLEXES}

\centerline{\bf CHIRAL JANUS COMPLEXES}

\

%\centerline{\it Mother of all acyclic complexes }

\

\centerline{Fyodor Malikov and Vadim Schechtman}

\bigskip\bigskip

\centerline{October 6, 2022}

\

\

\centerline{\bf Abstract}

\

\

%We make some elementary remarks concerning the relationships between chain complexes of infinite spheres and Koszul 
%complexes, and  chiral versions thereof. 

We propose chiral analogues of some infinite complexes appearing in the description of the coherent 
derived categories for projective spaces.

\

\

\centerline{\bf Introduction}

\

\

The aim of this note is to propose some examples of acyclic complexes unbounded in both directions.
 
One might ask why do we care about acyclic complexes, after all they are ignored when we define derived 
categories? One possible motivation is a remarkable but not very wellknown description, due to Beilinson and Bernstein, of the coherent derived category $\CD(\BP^n)$  over a projective space, given in [G]. Namely, 
the last category turns out to be equivalent to a category $\CA_n$ whose objects are unbounded acyclic complexes 
of free modules over the Grassmann (exterior) algebra $\Lambda[\BC^{n+1}]$, cf. {\it op. cit.} \S 5. Here it is important that  
the complexes are unbounded in both directions.

Examples which we discuss in \S \S 1 - 3 of this note  concern the simplest case $n = 0$.  So $\BP^0$ is just a point, and   
coherent sheaves over it are finite dimensional vector spaces. We consider a complex $\CK os\in \CA_0$ 
representing a one-dimensional vector space. The Grassmann algebra $\Lambda[\BC]$ is the algebra of dual numbers $\BC[\xi], \xi^2 = 0$.

The complex $\CK os$ which we call a {\it Janus} Koszul complex (cf. [P], [S], [KS] for a discussion of Janus objects) is defined as follows in {\it loc. cit.} (In the main body of the paper we will use,  
following {\it loc. cit.}, the notation $\Delta$ for this complex.)  
 
Its right part $\CK os_+$ is concentrated in degrees $\geq 0$ and is by definition the Koszul 
resolution for the $\BC[x]$-module $\BC$ where we assign to $x$ the degree $1$. As a vector space
$$
\CK os_+ = \BC[\xi]\otimes \BC[x],
$$
and its $i$-th term is a two-dimensional space $\BC[\xi]x^i$. It has only one nonzero cohomology, 
in degree $0$, which is $1$-dimensional.

The left part, 
$\CK os_-$, is concentrated in degrees $\leq - 1$ and is a sort of dual to $\CK os_-$; it has 
a unique nonzero cohomology, in degree $-1$, which is $1$-dimensional as well. 
Afterwards these two guys are docked into one complex, like Apollo - Soyuz.  
 
%Alternatively $\CK os$ may be defined in one shot: namely by replacing   
%the polynomial ring $\BC[x]$ by Laurent polynomials $\BC[x,x^{-1}]$ in the above 
%definition of $\CK os_+$.
 
The whole complex admits a {\it d\'ecalage} symmetry
$$
\sigma: \CK os\iso \CK os[1]. 
$$
For details see \S 2.

We remark that $\CK os_+$ may be identified with the cochain complex of a cellular complex 
homeomorphic to the infinite sphere $S^\infty$. Such complexes are discussed in \S 1, which is accessible to 
high school students.

In \S 3 we describe a chiral analogue $^\chr\CK os$ of $\CK os$.   

The chiral counterpart of $\CK os_+$ 
is a dg vertex algebra $^\chr\CK os_+$, and the chiral counterpart $^\chr\CK os_-$ of $\CK os_-$ is  vertex contragradient dual module for $^\chr\CK os_+$, see 3.5 and [FBZ], 10.4.6. The whole $^\chr\CK os$ is 
obtained by docking the two halves, similarly to the finite dimensional case.

%It is a vertex algebra\footnote{$^\chr\CK os$ is not a vertex algebra. Only half of it is, over which another is a module; so you can say it is a vertex algebra module, but even then the differential on the ``bouquet'' does not respect this structure -- very much like the finite dimensional counterpart. I suggest we simply remove this sentence.}  equipped with a differential, nonnegatively graded by the conformal weight, whose conformal weight zero component is $\CK os$. 

In \S 4 the above construction is generalized to the case $n > 0$. In addition, we present in 4.3 an interpretation of  
the negative part of the chiral Janus Koszul as a local cohomology of its positive part.

\

We are grateful to A.Beilinson for discussions and bringing [G] to our attention, and to M.Kapranov 
for interesting comments and an important correction.  

%\centerline{Quelques exercices}

%\centerline{\bf Complexe de de Rham}

\begin{CJK}{UTF8}{min}

%\centerline{ヘッケの縮退代数とコ骸骨\footnote{hecke no shukutai daisu to kogaikotsu}}
%\centerline{ニルヘッケの代数とコ骸骨\footnote{niruhecke no daisu to kogaikotsu}}

\end{CJK}

\bigskip\bigskip

\newpage

%\centerline{\bf \S 1.  AN YNFINITE SPHERE}

\centerline{\bf \S 1. An infinite sphere}

\

\

{\bf 1.1. Spheres.} We fix a base commutative ring $\fl$.

{\bf 1.1.1. The circle $S^1$}. Consider a cellular decomposition of 
$$
S^1 = \{ z\in \BC|\ |z| = 1\}
$$
having two $0$-cells
$$
c_0^{+} = 1,\ c_0^{-} = - 1.
$$
and two $1$-cells 
$$
c_1^{+} = \{e^{i\theta}| \ 0\leq \theta \leq \pi\}, \ 
c_1^{-} = \{e^{i\theta}| \ \pi\leq \theta \leq 2\pi\}, 
$$
where $c_1^{+}$ is 
oriented clockwise, and $c_1^{-}$ - counterclockwise. 

%as shewn below.

%\bigskip\bigskip

%\includegraphics[width = 10cm]{cercle.jpg}

%\centerline{Fig. 1. The circle.}

\

So we have defined a $CW$-complex, to be denoted $K_{[0,1]}$, whose complex of chains is 
$$
C_\bullet(K_{[0,1]}; \fl):\ 0 \lra C_1(K_{[0,1]}; \fl)\overset{d}\lra C_0(K_{[0,1]}; \fl)\lra 0
$$
where
$$
C_1(K_{[0,1]}; \fl) = \fl c_1^+\oplus \fl c_1^-,\ C_0(K_{[0,1]}; \fl) = \fl c_0^+\oplus \fl c_0^-,
$$
and the differential 
$$
dc_1^+ = dc_1^- = c_0^+ - c_0^-,
$$
i.e. $d$ is given by matrix   
$$
A = \left(\begin{matrix} 1 & -1\\
1 & -1
\end{matrix}\right).
\eqno{(1.1.1)}
$$
We define a complex $\CC^\bullet_{[0,1]}$ as $C_\bullet(K_{[0,1]}; \fl)$ placed in degrees $0, 1$.

Thus 
$$
H^0(\CC^\bullet_{[0,1]})\isom\fl;\ H^1(\CC^\bullet_{[0,1]})\isom\fl
$$

\

{\bf 1.1.2. The sphere $S^2$.} Similarly we define a cellular decomposition $K_{[0,2]}$ of the sphere $S^2$ 
having two cells $c_i^\pm$ in dimensions $i = 0, 1, 2$. Namely the $2$-cells $c_2^\pm$ are northern and southern  hemispheres, and 
$$
K_{[0,1]}\subset K_{[0,2]}
$$
as the equator. Its complex of chains will be of length $2$ 
$$
C_\bullet(K_{[0,2]}; \fl):\ 0 \lra C_2(K_{[0,2]};\fl)\overset{d_2}\lra C_1(K_{[0,2]};\fl)\overset{d_1}\lra C_0(K_{[0,2]};\fl)\lra 0
$$
with
$$
C_i(K_{[0,2]};\fl) = \fl c_i^+\oplus \fl c_i^-,
$$
and both differentials are given by the matrix $A$ (we note that $A^2 = 0$). 

We define a complex $\CC^\bullet_{[0,2]}$ as $C_\bullet(K_{[0,2]}; \fl)$ placed in degrees $0, 1, 2$.

Thus 
$$
H^0(\CC^\bullet_{[0,2]})\isom\fl,\ H^1(\CC^\bullet_{[0,2]}) = 0,\ H^2(\CC^\bullet_{[0,2]})\isom\fl.
$$

\

{\bf 1.1.3. The sphere $S^n$.} Similarly for an arbitrary $n$ we define a cellular decomposition 
$K_{[0,n]}$ of the sphere $S^n$ 
having two cells $c_i^\pm$ in dimensions $i = 0, 1, \ldots, n$, with the  
 embedding
$$
K_{[0,n-1]}\subset K_{[0,n]}
$$
as the equator. 

Its complex of chains will be of length $n$ 
$$
C_\bullet(K_{[0,n]}; \fl):\ 0 \lra C_n(K_{[0,n]};\fl)\overset{d_n}\lra \ldots \overset{d_1}\lra C_0(K_{[0,n]};\fl)\lra 0
$$
with
$$
C_i(K_{[0,n]};\fl) = \fl c_i^+\oplus \fl c_i^-,
$$
and both differentials given by the matrix $A$. 

We define a complex $\CC^\bullet_{[0,n]}$ as $C_\bullet(K_{[0,n]}; \fl)$ placed in degrees $0, \ldots, n$.

Thus 
$$
H^0(\CC^\bullet_{[0,n]})\isom\fl,\ H^i(\CC^\bullet_{[0,2]}) = 0, 1\leq i \leq n-1,\ 
H^n(\CC^\bullet_{[0,n]})\isom\fl.
$$

\

{\bf 1.2.} As remarked M.Kapranov, we may consider $K_{[0,n]}$ as a {\it globular $n$-category} 
(cf. for example [ManS]) having two objects, two $1$-morphisms, two $2$-morphisms, etc. 

There exists a canonical functor of $n$-categories 
$$
\phi_n:\ K_{[0,n]} \lra HB_{n+2}
$$
where $HB_{n+2}$ is the globular $n$-category of higher Bruhat orders on the symmetric group $\Sigma_{n+2}$.
This follows from the main result of {\it op. cit.}

For example $\phi_1$ takes two objects of $K_{[0,1]}$ to permutations $e= ()$ and $(31)$ from $\Sigma_3$, and two arrows 
to two possible paths from $e$ to $(13)$ for the weak Bruhat order in $\Sigma_3$.   

\

{\bf 1.3. The space $S^\infty$.} We may pass to the limit $n\lra\infty$, and define an infinite dimensional 
$CW$-complex
$$
K_{[0,\infty[} = \cup_{n=1}^\infty K_{[0,n]}.
$$
It is contractible.

The corresponding complex 
$$
\CC_+^\bullet = \CC^\bullet_{[0,\infty[}
$$
will live in degrees $\geq 0$: namely 
$$
\CC_+^i = \fl c_+^i\oplus \fl c_-^i,\ i\geq 0,\ 
$$
with $d = A$, 
and will have the cohomology
$$
H^0(\CC^\bullet_{[0,\infty[}) \isom\fl,\ 
H^i(\CC^\bullet_{[0,\infty[}) = 0,\ i > 0.
$$
Similarly, using the left shifts of complexes $C_\bullet(K_{[0,n]};\fl)$ we define a complex 
$$
\CC_-^\bullet = \CC^\bullet_{]-\infty, 0]}
$$
living in degrees $\leq 0$ with 
$$
\CC_-^i = \fl c_+^i\oplus \fl c_-^i,\ i\leq 0,\ 
$$
with $d = A$, 
having  the cohomology
$$
H^0(\CC^\bullet_{]-\infty, 0]}) \isom\fl,\ 
H^i(\CC^\bullet_{]-\infty, 0]}) = 0,\ i < 0.
$$  

\

{\bf 1.4. Double infinite, or Janus.} We define a complex 
$$
\CC^\bullet = \CC^\bullet_{]-\infty, \infty[}
$$
living in all degrees $i\in \BZ$ with 
$$
\CC^i = \fl c_+^i\oplus \fl c_-^i,\ i\in \BZ,\ 
$$
with $d = A$. 

It is acyclic:
$$
H^i(\CC^\bullet) = 0,\ i\in \BZ.
$$

\

\newpage

\centerline{\bf \S 2. Koszul complex}

\

{\bf 2.1.} The following is a particular case of a construction due to A.Beilinson and J.Bernstein, see [G], \S 5.  

Let 
$\Xi$ be a one-dimensional vector space over $\fl$ with  base $\xi$, $X = \Xi^*$ the dual space 
with the base $x$ dual to $\xi$. Let
$$
\Lambda(\Xi) = \fl\oplus \Xi
$$
be the exterior algebra, $\Lambda(\Xi)[-1]$ the free $\Lambda(\Xi)$-module of rank $1$ with a generator 
$e$ living in degree $-1$, 
$$
S^\bullet(X) = \oplus_{i=0}^\infty \fl x^i
$$
the symmetric algebra. 

{\bf 2.1.1. Positive part.} Let 
$$
\Delta_+^\bullet = \Lambda(\Xi)[-1]\otimes S^\bullet(X) 
$$
be the shifted Koszul complex. By definition it is a complex of $\Lambda(\Xi)$-modules 
$$
\Delta_+^\bullet:\ 0\lra \Lambda(\Xi)[-1]\otimes S^0(X)\overset{d_+}\lra \Lambda(\Xi)[-1]\otimes S^1(X)
\overset{d_+}\lra 
\ldots
$$
with a differential
$$
d(\mu\otimes t) = \xi\mu\otimes xt
$$
The only nonzero cohomology of this complex is 
$$
H^0(\Delta_+^\bullet)\isom \fl.
$$ 

\ 

{\bf 2.1.2. Negative part.} For a $\Lambda(\Xi)$-module $M$ denote  
$$
M^* = Hom_{\Lambda(\Xi)}(M, \Lambda(\Xi))
$$
the dual module. 
The complex $\Delta_-^\bullet$ consists of dual $\Lambda(\Xi)$-modules and lives in degrees 
$\leq - 1$:
$$
\Delta_-^\bullet: \ldots \overset{d_-}\lra \Delta_-^{-2}\overset{d_-}\lra \Delta_-^{-1} \lra 0 
$$
where 
$$
\Delta_-^{-i} = (\Delta_+^{i-1})^* = (\Lambda(\Xi)[-1])^*\otimes S^i(X)
$$
and
$$
d_-^{-i} = (d_+^{i-2})^*
$$
The only nonzero cohomology of this complex is 
$$
H^{-1}(\Delta_-^\bullet)\isom \fl.
$$ 

\

{\bf 2.1.3. Janus complex $\Delta^\bullet$.} We glue the above two guys together using a map of 
$\Lambda(\Xi)$-modules
$$
d_-^{-1}:\ \Delta_-^{-1} = (\Lambda(\Xi)[-1])^* = \Lambda(\Xi)\lra \Lambda(\Xi)[-1] = 
\Delta_+^0,
$$
$$
d_-^{-1}(1) = \xi,
$$ 
and get an infinite in both sides acyclic complex
$$
\Delta^\bullet: \ldots \overset{d_-}\lra \Delta_-^{-2}\overset{d_-}\lra \Delta_-^{-1} 
\lra \Delta_+^0 \overset{d_+}\lra \Delta_+^{1} \lra \ldots    
$$

\

{\bf 2.2. Localized Koszul.}  The Koszul of 2.1.1, $\Delta_+^\bullet$, is often realized slightly differently. Given an $\fl$-algebra $R$ with an element $f\in R$,  the Koszul complex can be defined to be the algebra  $\fl[\xi]\otimes R$ with differential defined to be the following derivation
\[
d_K=\partial_\xi \otimes f,
\]
where
$$
\dpar_\xi(\xi) = 1,\ \dpar_\xi(1) = 0.
$$

Denote this object by $K(R,f)$.

If $R=\fl[x]$, then we give it a grading by letting $K(\fl[x],f)^i=\fl[\xi]x^i$.

It is easy to see that $K(\fl[x],x)$ is isomorphic with $\Delta_+^\bullet$ as a complex, the isomorphism 
\[
K(\fl[x],x)\rightarrow  \Delta_+^\bullet
\]
being defined by the assignment
\[
1\otimes x^n\mapsto\xi\otimes x^n;\; \xi\otimes x^n\mapsto 1\otimes x^n.
\]

In the case where $R=\fl[x,x^{-1}]$ (the grading being similarly defined by the powers of $x$), the Koszul complex
$K(\fl[x,x^{-1}],x)$ is obviously acyclic and is  in fact  naturally isomorphic, as a complex, with $\Delta^\bullet$,  the  Janus complex of 2.1.3.   To define this isomorphism choose a basis of $\fl[\xi]\otimes\fl[x]$ to be the set of monomials $\{1\otimes x^n,\xi\otimes x^n,\; n\geq 0\}$ and consider a map 
\[
\fl[\xi]\otimes\fl[x^{-1}]x^{-1}\rightarrow (\fl[\xi]\otimes\fl[x])^*
\]
defined by the assignment:
\[
1\otimes x^{-n-1}\mapsto (\xi\otimes x^n)^\vee,\; \xi\otimes x^{-n-1}\mapsto (1\otimes x^n)^\vee,
\]
where $\{(1\otimes x^n)^\vee,(\xi\otimes x^n)^\vee,\; n\geq 0\}$ is the basis dual to $\{1\otimes x^n,\xi\otimes x^n,\; n\geq 0\}$.

Using the inverse dual to the  above defined isomorphism $K(\fl[x],x)\rightarrow  \Delta_+^\bullet$ we obtain the composite map
\[
\fl[\xi]\otimes\fl[x^{-1}]x^{-1}\rightarrow (\fl[\xi]\otimes\fl[x])^*\rightarrow  (\Delta_+^\bullet)^*=\Delta_-^\bullet,
\]
cf. 2.1.2, which is easily seen to be an isomorphism of complexes in degrees -2,-3,-4, etc.  Taking the direct sum we obtain a vector
space isomorphism
\[
K(\fl[x,x^{-1}],x)=\fl[\xi]\otimes\fl[x^{-1}]x^{-1}\oplus \fl[\xi]\otimes\fl[x]\rightarrow \Delta^\bullet_-\oplus \Delta^\bullet+=\Delta^\bullet.
\]

Under this isomorphism the gluing map, $\Delta^{-1}_-\rightarrow\Delta^0_+$, which in 2.1.3 had to be defined by hand,  is on the left  simply the restriction of  the standard Koszul differential to $\fl[\xi]\otimes x^{-1}$:
\[
\partial_\xi\otimes x(\xi\otimes x^{-1})=1\otimes 1,\; \partial_\xi\otimes x(1\otimes x^{-1})=0.
\]

Note that, unlike $\Delta_+^\bullet$ or  $\Delta^\bullet$, $K(R,f)$, $K(\fl[x],x)$, $K(\fl[x,x^{-1}],x)$ are   commutative dg-algebras.

\

{\bf 2.3. Comparison.}
The matrix of $d_+$ is 
$$
B = \left(\begin{matrix} 0 & 0\\
1 & 0
\end{matrix}\right).
%\eqno{(2.2.1)}
$$
The matrices $A$ and $B$ are conjugated: if
$$
\alpha = \left(\begin{matrix} 1\\ 0
\end{matrix}\right),\ \beta = \left(\begin{matrix} 1\\ 1
\end{matrix}\right)
$$
then
$$
A\alpha = \beta,\ A\beta = 0.
$$
Thus a base change gives rise to an isomorphism of complexes 
$$
\Delta_+^\bullet\isom \CC^\bullet_+
$$
Similarly
$$
\Delta_-^\bullet\isom \CC^\bullet_-,
$$
and
$$
\Delta^\bullet\isom \CC^\bullet.
$$

\

\

%\centerline{\bf \S 3. CHIRAL KOSZUL OR DE RHAM COMPLEX}

\centerline{\bf \S 3. Chiral Koszul: one variable }

%\centerline{\bf one variable}

\

{\bf 3.1.}  From now on our ground ring $\fl$ will be $\BC$ (or any commutative ring containing 
$\BQ$).  

Consider $T^*\BC$ as an algebraic variety. Its ring of regular functions, $\BC[T^*\BC]$, is a polynomial ring on 2 variables,
$\BC[x,\partial_x]$, where $x$ is a coordinate on $\BC$ and $\partial_x$ is the vector field s.t. $\partial_x(x)=1$. In fact,
$\BC[T^*\BC]$ is naturally a Poisson algebra with Poisson bracket bracket uniquely determined by the condition $\{\partial_x,x\}=1$,
$\{\partial,\partial\}=\{x,x\}=0$.

Now consider the jet-space, $J_\infty T^*\BC$. Its ring of regular functions $\BC[J_\infty T^*\BC]$ is a polynomial ring on infinitely many variables, $\BC[x_{-(n)},\partial_{x,(-n)},\,n>0]$, with derivation, which is uniquely determined by its action on the generators,
$x_{(-n)}\mapsto nx_{(-n-1)}$, $\partial_{x,(-n)}\mapsto n\partial_{x,(-n-1)}$.

$\BC[J_\infty T^*\BC]$ is not a Poisson algebra in any natural way; it is instead a {\em coisson} or vertex Poisson algebra. The coisson,
or Lie*,  bracket is determined by the assignment
\[
\{\partial_{x,(-n)},x_{(-m)}\}=\frac{1}{(n-1)!(m-1)!}\partial_z^{n-1}\partial_w^{m-1}\delta(z-w),
\]
where we use some standard notation without explanation, referring the reader  to the inspired and inspirational Introduction to [BD].

This coisson algebra structure is easy to quantize and the result of this quantization is well known.  As a vector space, it is still the ring
 $\BC[x_{-(n)},\partial_{x,(-n)},\,n>0]$ and the state-field correspondence is determined by the assignment
\begin{eqnarray*}
1&\mapsto& \mbox{Id},\\
x_{(-n)}&\mapsto&\frac{1}{(n-1)!}\partial_z^{n-1}x(z),\\
\partial_{x,(-n)}&\mapsto&\frac{1}{(n-1)!}\partial_z^{n-1}\partial_x(z),
\end{eqnarray*}
where the only nontrivial OPE among the fields $\partial_x(z),x(x)$ is as follows:
\[
\partial_x(z)x(w)=\frac{1}{z-w}+\cdots
\]
This vertex algebra is often referred to as the ``$\beta\gamma$-system.'' We  called it an algebra of chiral differential operators (CDO)
and denoted it $D_\BC^{ch}$, cf. [GMS]. 

\begin{sloppypar}{\bf 3.2.} This construction has an obvious superalgebra extension. The underlying jet-space is $J_\infty T^*\Pi T\BC$. The corresponding ring of functions is the supercommutative algebra $\BC[x_{(-n)},\partial_{x,(-n)}, \xi_{(-n)},\partial_{\xi,(-n)},\,n>0]$ with derivation
$x_{(-n)}\mapsto nx_{(-n-1)}$, $\partial_{x,(-n)}\mapsto n\partial_{x,(-n-1)}$, $\xi_{(-n)}\mapsto n\xi_{(-n-1)}$, $\partial_{\xi,(-n)}\mapsto n\partial_{\xi,(-n-1)}$; the variables $\xi_{(-n)}$ and $\partial_{\xi,(-n)}$ are odd.
\end{sloppypar}

This  is a supercoisson algebra, and its quantization is similar to the one described in sect. 3.1: the space is  $\BC[x_{(-n)},\partial_{x,(-n)}, \xi_{(-n)},\partial_{\xi,(-n)},\,n>0]$, and the state-field correspondence is obtained by extending the assignment at the end  of 3.1 as follows:
\begin{eqnarray*}
\xi_{(-n)}&\mapsto&\frac{1}{(n-1)!}\partial_z^{n-1}\xi(z),\\
\partial_{\xi,(-n)}&\mapsto&\frac{1}{(n-1)!}\partial_z^{n-1}\partial_\xi(z),
\end{eqnarray*}
The fields $\partial_\xi(z)$, $\xi(z)$ are odd, and the only nontrivial OPE they satisfy is this:
\[
\partial_\xi(z)\xi(w)=\frac{1}{z-w}+\cdots
\]
This is a super CDO, to be denoted $D^{ch}_{\BC^{1|1}}$.

{\bf 3.3.}  One advantage of dealing with a super CDO is that it contains various differentials. For example, consider the map
\[
d^{ch}_{K}: D^{ch}_{\BC^{1|1}}\rightarrow D^{ch}_{\BC^{1|1}}
\]
defined by 
\[
d^{ch}_{K}=\int\partial_{\xi}(z)x(z)\,dz.
\]
If we let the degree of $x_{(n)}$ be one, of $\partial_{x,(n)}$ be -1, and that of $\xi_{(n)}$ and $\partial_{\xi,(n)}$ be 0, for any $n$, then one easily verifies that $d^{ch}_{K}$ is of square 0 and  degree 1, and so the pair $(D^{ch}_{\BC^{1|1}},d^{ch}_{K})$ is a complex infinite in both
directions -- in fact a {\em differential graded vertex algebra} as the integral of a field, such as $d^{ch}_{K}$, is well known to be a derivative of a vertex algebra structure. One also notices that 
\[
\BC[x_{(-1)},\xi_{(-1)}]= \BC[x_{(-1)}]\oplus \BC[x_{(-1)}]\xi_{(-1)}\subset D^{ch}_{\BC^{1|1}}
\]
 is a subcomplex, and as such it is rather naturally identified with the Koszul complex $K(\fl[x],x)$ of 2.2  if we think of $x_{(-1)}$ as the
 coordinate $x$ on $\BC$, and of $\xi_{(-1)}$ as $\xi$. For example, the restriction of $d^{ch}_K$ to $\BC[x_{(-1)},\xi_{(-1)}]$ is this:
 \[
 d^{ch}_K|_{\BC[x_{(-1)},\xi_{(-1)}]}=x_{(-1)}\frac{d}{d\xi_{(-1)}},
 \]
 which of course coincides with $d_K$ of 2.2.

 Furthermore, one notices that thus defined  embedding $K(\fl[x],x)\hookrightarrow D^{ch}_{\BC^{1|1}}$ is a quasiisomorphism. This is because of a part of an N=2 superconformal structure that $D^{ch}_{\BC^{1|1}}$ carries.  Consider 2 more fields
 \[
 L(z)=(x_{(-2)}\partial_{x,(-1)})(z)+ (\xi_{(-2)}\partial_{\xi,(-1)})(z),\;
 G(z)=(\xi_{(-2)}\partial_{x,(-1)})(z).
 \]
 One readily verifies that $L(z)$ is the Virasoro field and, in particular, the component $L_{(1)}=\int L(z)z\,dz$, as an operator acting on
 $D^{ch}_{\BC^{1|1}}$, defines the conformal grading,
 \[
 D^{ch}_{\BC^{1|1}}=\bigoplus_{n=0}^\infty D^{ch}_{\BC^{1|1},n},
 \]
 where $D^{ch}_{\BC^{1|1},n}$ is the eigenspace of $L_{(1)}$ of eigenvalue $n$. A routine and familiar computation will show that this grading is, in fact, a grading on the underlying polynomial ring,
  \[
 D^{ch}_{\BC^{1|1}}= \BC[x_{(-n)},\partial_{x,(-n)}, \xi_{(-n)},\partial_{\xi,(-n)},\,n>0],
 \]
 where 
 \[
 \mbox{deg}x_{(-n)}=n-1;\; \mbox{deg}\partial_{x,(-n)}=n;\;  \mbox{deg}\xi_{(-n)}=n-1;\; \mbox{deg}\partial_{\xi,(-n)}=n.
 \]

 This grading is preserved by the differential $d^{ch}_K$.
 For example, the conformal weight zero subspace, $D^{ch}_{\BC^{1|1},0}$ is exactly the ordinary Koszul complex $K(\fl[x],x)$ in the form of
 $\BC[x_{(-1)},\xi_{(-1)}]$,  which we considered above, $x_{(-1)},\xi_{(-1)}$, being the only generators of conformal weight 0, and the differential $d^{ch}_{K}$ preserves this grading.
 
 Furthermore,
 \[
 [d^{ch}_{K}, \int G(z)z\,dz]=L_{(1)},
 \]
 which means that the restriction of the differential $d^{ch}_{K}$ to each {\em nonzero} conformal weight subspace is homotopic to identity. This implies that  the embedding $K(\fl[x],x)\hookrightarrow D^{ch}_{\BC^{1|1}}$ is indeed a quasiisomorphism.
 
 Therefore, the cohomology $H_{d^{ch}_{K}}(D^{ch}_{\BC^{1|1}})$ is 1-dimensional and is spanned by the class
 of $1\in D^{ch}_{\BC^{1|1}}$.
 
 \bigskip
 
 {\bf 3.4.}  We can localize and consider 
 \[
 D^{ch}_{\BC^{1|1}}[x_{(-1)}^{-1}]\stackrel{\mbox{def}}{=} \BC[x_{(-1)}, x_{(-1)}^{-1}, x_{(-n-1)},\partial_{x,(-n)}, \xi_{(-n)},\partial_{\xi,(-n)},\,n>0].
 \]
 It is still a vertex algebra, [MSV,GMS],  and it shares many features with the above analyzed $ D^{ch}_{\BC^{1|1}}$:
 
 (1) the pair $(D^{ch}_{\BC^{1|1}}[x^{-1}_{(-1)}],d^{ch}_K)$ is a differential graded vertex algebra with differential $d^{ch}_K$;
 
 (2) it is conformally graded by the eigenvalues of $L_{(1)}$:
  \[
 D^{ch}_{\BC^{1|1}}[x_{(-1)}^{-1}]=\bigoplus_{n=0}^\infty D^{ch}_{\BC^{1|1},n}[x_{(-1)}^{-1}],
 \]
 the grading being determined in an obvious manner by that of the polynomial ring $ \BC[x_{(-n)},\partial_{x,(-n)}, \xi_{(-n)},\partial_{\xi,(-n)},\,n>0]$, see 3.3;
 
 (3) localization of the identification $K(\fl[x],x)= D^{ch}_{\BC^{1|1},0}$, see 3.3, gives an identification 
 $K(\fl[x,x^{-1}],x)= D^{ch}_{\BC^{1|1},0}[x_{(-1)}^{-1}]$; and, finally,
 
 (4) the restriction of $d^{ch}_K$ to $D^{ch}_{\BC^{1|1},n}[x_{(-1)}^{-1}]$ with $n>0$ is homotopic to identity (thanks to the relation $ [d^{ch}_{K}, \int G(z)z\,dz]=L_{(1)}$) and, therefore, the emebedding  (as the conformal weight zero subspace, see point (3))
 \[
 K(\fl[x,x^{-1}],x)\hookrightarrow D^{ch}_{\BC^{1|1}}[x_{(-1)}^{-1}]
 \]
 is a quasiisomorphism.
 
 It is in this sense that we assert $ D^{ch,}_{\BC^{1|1}}[x_{(-1)}^{-1}]$ is a chiralization of $ K(\fl[x,x^{-1}],x)$, hence of the 
 Janus $\Delta^\bullet$, 2.1.3.

 {\bf 3.5.   The chiral  Janus complex. } An obvious embedding $D^{ch}_{\BC^{1|1}}\hookrightarrow D^{ch}_{\BC^{1|1}}[x_{(-1)}^{-1}]$
 is a vertex algebra morphism; therefore, the quotient $D^{ch}_{\BC^{1|1}}[x_{(-1)}^{-1}]/D^{ch}_{\BC^{1|1}}$ is a $D^{ch}_{\BC^{1|1}}$-module - in fact a diffferential graded $D^{ch}_{\BC^{1|1}}$-module as $d^{ch}_K$ tautologically acts on any $D^{ch}_{\BC^{1|1}}$-module. 
 
 $D^{ch}_{\BC^{1|1}}[x_{(-1)}^{-1}]/D^{ch}_{\BC^{1|1}}$ is conformally graded by the eigenvalues of $L_{(1)}$, because $D^{ch}_{\BC^{1|1}}[x_{(-1)}^{-1}]$ and $D^{ch}_{\BC^{1|1}}$ are,  3.3, 3.4, and the grading is preserved by $d^{ch}_K$, which implies, as in 3.3, that
  the cohomology $H_{d^{ch}_K}(D^{ch}_{\BC^{1|1}}[x_{(-1)}^{-1}]/D^{ch}_{\BC^{1|1}})$ is 1-dimensional, too, and  is spanned by the cocycle $\xi_{(-1)}/x_{(-1)}$, just as $H_{d^{ch}_K}(D^{ch}_{\BC^{1|1}}[x_{(-1)}^{-1}]/D^{ch}_{\BC^{1|1}})$ is spanned by the cocycle $1$, the last sentence of 3.3.
 
 We can now mimic 2.1.3 and form a bouquet of the 2 complexes, $D^{ch}_{\BC^{1|1}}[x_{(-1)}^{-1}]/D^{ch}_{\BC^{1|1}}$ and $D^{ch}_{\BC^{1|1}}$  as follows. Define 
 \[
 D^{ch}_{\BC^{1|1}}[x_{(-1)}^{-1}]/D^{ch}_{\BC^{1|1}}\stackrel{.}{\oplus}D^{ch}_{\BC^{1|1}}
 \]
to be
\[
D^{ch}_{\BC^{1|1}}[x_{(-1)}^{-1}]/D^{ch}_{\BC^{1|1}}\oplus D^{ch}_{\BC^{1|1}}
\]
as a vector space. To define the differential we, first, stipulate  that in nonzero conformal weight subspaces  it is the direct
sum of the original differentials in nonzero conformal weight subspaces. More formally, on 
\[
D^{ch}_{\BC^{1|1},n}[x_{(-1)}^{-1}]/D^{ch}_{\BC^{1|1,0}}\oplus D^{ch}_{\BC^{1|1},n},\; n>0,
\]
the differential is
\[
d^{ch}_K\oplus d^{ch}_K:\; 
D^{ch}_{\BC^{1|1},n}[x_{(-1)}^{-1}]/D^{ch}_{\BC^{1|1,0}}\oplus D^{ch}_{\BC^{1|1,n},0}\rightarrow
D^{ch}_{\BC^{1|1},n}[x_{(-1)}^{-1}]/D^{ch}_{\BC^{1|1,0}}\oplus D^{ch}_{\BC^{1|1,n},0} .
\]

Now focus on the 0 conformal weight subspace,
\[
D^{ch}_{\BC^{1|1},0}[x_{(-1)}^{-1}]/D^{ch}_{\BC^{1|1},0}\oplus D^{ch}_{\BC^{1|1},0}.
\]
By construction, $D^{ch}_{\BC^{1|1},0}$ is the ordinary Koszul, $K(\fl[x],x)$, see an explanation of this in 3.3, and 
$D^{ch}_{\BC^{1|1},0}[x_{(-1)}^{-1}]/D^{ch}_{\BC^{1|1,0}}$ is  its dual, see 2.2. This places us entirely in the situation
of 2.2,  where, in particular $D^{ch}_{\BC^{1|1},0}$ is graded by positive powers of $x_{(-1)}$,
\[
D^{ch}_{\BC^{1|1},0}=D^{ch}_{\BC^{1|1},0}[0]\oplus D^{ch}_{\BC^{1|1},0}[1]\oplus D^{ch}_{\BC^{1|1},0}[2]\oplus\cdots,
\]
$D^{ch}_{\BC^{1|1},0}[x_{(-1)}^{-1}]/D^{ch}_{\BC^{1|1,0}}$  by negative,
\[
D^{ch}_{\BC^{1|1},0}[x_{(-1)}^{-1}]/D^{ch}_{\BC^{1|1,0}}=
\cdots \oplus (D^{ch}_{\BC^{1|1},0}[x_{(-1)}^{-1}]/D^{ch}_{\BC^{1|1,0}})[-2]\oplus
(D^{ch}_{\BC^{1|1},0}[x_{(-1)}^{-1}]/D^{ch}_{\BC^{1|1,0}})[-1].
\]
In either case, the $n$-th graded component is 2-dimensional with basis $\{\xi_{(-1)}x_{(-1)}^n, x_{(-1)}^n\}$.

We now mimic 2.1.3, as we did in 2.2, to join these complexes by using the map
\[
(D^{ch}_{\BC^{1|1},0}[x_{(-1)}^{-1}]/D^{ch}_{\BC^{1|1,0}})[-1]\rightarrow D^{ch}_{\BC^{1|1},0}[0],
\]
which sends
\begin{eqnarray*}
(D^{ch}_{\BC^{1|1},0}[x_{(-1)}^{-1}]/D^{ch}_{\BC^{1|1,0}})[-1]\ni \frac{\xi_{(-1)}}{x_{(-1)}}&\mapsto& 1\in D^{ch}_{\BC^{1|1},0}[0],\\
(D^{ch}_{\BC^{1|1},0}[x_{(-1)}^{-1}]/D^{ch}_{\BC^{1|1,0}})[-1]\ni \frac{1}{x_{(-1)}}&\mapsto& 0\in D^{ch}_{\BC^{1|1},0}[0].
\end{eqnarray*}
The result is the bouquet $D^{ch}_{\BC^{1|1}}[x_{(-1)}^{-1}]/D^{ch}_{\BC^{1|1}}\stackrel{.}{\oplus}D^{ch}_{\BC^{1|1}}$. It is acyclic by construction.

\bigskip

{\bf 3.6. Remark.}   $D^{ch}_{\BC^{1|1}}$ is a conformal vertex algebra; as such it defines a sheaf on $\BP^1$. In this situation, one defines a duality functor on the category of $ D^{ch}_{\BC^{1|1}}$-modules, see [FBZ], 10.4.6. $D^{ch}_{\BC^{1|1}}[x_{(-1)}^{-1}]/D^{ch}_{\BC^{1|1}}$ is naturally a $ D^{ch}_{\BC^{1|1}}$-module and, in fact, it is the dual to $ D^{ch}_{\BC^{1|1}}$.

\bigskip

{\bf 3.7.}  Let us explore the chiral Janus complex $D^{ch}_{\BC^{1|1}}[x_{(-1)}^{-1}]/D^{ch}_{\BC^{1|1}}\stackrel{.}{\oplus}D^{ch}_{\BC^{1|1}}$ in some detail.  Its one essential feature is that although infinite in both directions, as a complex, it is a direct sum of conformal weight subspaces,  3.3, 3.4, and each such subspace is a subcomplex bounded in both directions.We shall now examine these subcomplexes, starting with the subcomplex  $D^{ch}_{\BC^{1|1}}$. 

It is convenient  to shift the subindices of  the generators in order to makes sure they respect the conformal grading. Therefore we re-denote
\[
x_{-n+1}=x_{(-n)},\partial_{x,-n}=\partial_{x,(-n)}, \xi_{-n+1}=\xi_{(-n)},\partial_{\xi,-n}=\partial_{\xi,(-n)},\; n>0.
\]
A moment's thought will show that indeed the conformal weight of either of the generators 
$x_{-n},\partial_{x,-n}, \xi_{-n},\partial_{\xi,-n}$ is $n$.

{\bf 3.7.1.} Consider the conformal weight 0 subspace, $D^{ch}_{\BC^{1|1},0}$. It is equal to the Laurent superpolynomial
ring, $\BC[x_0,\xi_0]$, hence as a complex it is identified with the de Koszul complex $K(\BC[x_0], x_0)$, cf. 3.3.
 In fact, it is a direct sum
of 2-dimensional subcomplexes $C_n$, $n\geq 0$, $C_n$ being spanned by $\xi_0x_0^n, x_0^{n+1}$, and one 1-dimensional, $C_{-1}$, spanned by $1$. The matrix of the differential restricted to $C_n$, $n\geq 0$, in the indicated basis is
\[
\left(\begin{matrix} 0&0\\1&0\end{matrix}\right),
\]
cf. (2.2.1).  The element $1\in C_{-1}$ is a cocycle and spans the cohomology.

{\bf 3.7.2.} The conformal weight 1 component is a $\mbox{rk}=4$  free module  over the conformal weight 0 component, 
$\BC[x_0,\xi_0]$,
with basis $x_{-1},\partial_{x,-1},\xi_{-1},\partial_{\xi,-1}$.  Denote the 4-dimensional vector space spanned by these 4 vectors by $V_1$. We have 
\[
D^{ch}_{\BC^{1|1,1}}=V_1\otimes \BC[x_0,\xi_0],
\]
a tensor product of 2 complexes, the action of the differential on $V_1$ being given by
\begin{eqnarray*}
\xi_{-1} \mapsto x_{-1},\;x_{-1}\mapsto 0,\\
\partial_{x,-1} \mapsto -\partial_{\xi,-1},\;\partial_{\xi,-1} \mapsto 0.
\end{eqnarray*}
Of course this creates two more blocks of the familiar type
\[
\left(\begin{matrix} 0&0\\1&0\end{matrix}\right)
\]
along with yet another way of proving the acyclicity in conformal weight 1.

{\bf 3.7.3.} In an arbitrary conformal weight $N$ subspace, the situation is similar. We have a tensor product of complexes
\[
D^{ch}_{\BC^{1|1},N}=V_N\otimes \BC[x_0,\xi_0],
\]
where $V_N$ is a finite dimensional vector space with basis consisting of monomials $x_{-n_1}^b\partial_{x,-n_2}^c \xi_{-n_3}^d\partial_{\xi,-n_4}^e$, where the subindices $n_j$ are all strictly negative and satisfy
\[
bn_1+cn_2+dn_3+en_4=-N.
\]
The action of the differential is determined by the following conditions: 

$\bullet$ it is a derivation;

$\bullet$ it sends
\begin{align*}
\xi_{-n} \mapsto x_{-n},\;x_{-n}\mapsto 0,\\
\partial_{x,-n} \mapsto -\partial_{\xi,-n},\;\partial_{\xi,-n} \mapsto 0,
\end{align*}
which makes the acyclicity of $V_N$, hence that of $D^{ch}_{\BC^{1|1},N}$, $N>0$, manifest.

{\bf 3.7.4.} The subcomplex  $D^{ch}_{\BC^{1|1}}[x_{(-1)}^{-1}]/D^{ch}_{\BC^{1|1}}$ is dual to $D^{ch}_{\BC^{1|1}}$, which we have just
analyzed. In particular, its conformal weight 0 component is dual to  $D^{ch}_{\BC^{1|1},0}$, also has 1-dimensional cohomology, and the formation of the Janus  $D^{ch}_{\BC^{1|1},0}[x_{(-1)}^{-1}]/D^{ch}_{\BC^{1|1},0}\stackrel{.}{\oplus}D^{ch}_{\BC^{1|1},0}$  was discussed in 3.5. 

 Being isomorphic to $K(\BC[x_0],x_0)$, $D^{ch}_{\BC^{1|1},0}$ is also isomorphic to $\Delta^\bullet_+$, see 2.2, therefore 
 $D^{ch}_{\BC^{1|1},0}[x_{(-1)}^{-1}]/D^{ch}_{\BC^{1|1},0}\stackrel{.}{\oplus}D^{ch}_{\BC^{1|1},0}$ is nothing but the Janus complex of 2.1.3.

On the other hand, the higher conformal weight components of $D^{ch}_{\BC^{1|1}}[x_{(-1)}^{-1}]/D^{ch}_{\BC^{1|1}}$ are acyclic and enter the Janus
$D^{ch}_{\BC^{1|1}}[x_{(-1)}^{-1}]/D^{ch}_{\BC^{1|1}}\stackrel{.}{\oplus}D^{ch}_{\BC^{1|1}}$ simply as direct summands, 3.5.

\

\bigskip

\centerline{\bf \S 4. Chiral Koszul: several variables }

\bigskip

{\bf 4.1.} If $V_i$, $i\in I$, is a family of vertex algebras, then the tensor product over the ground ring, $\otimes_IV_i$, is canonically
a vertex algebra: if we write $v_i(z)$ for the field attached to $v_i\in V_i$, then we define the field attached to $\otimes_Iv_i$ to be
$\otimes_I(v_i(z))$. 

Just as the coordinate ring $\BC[\BC^I]=\otimes_I\BC[\BC]$ or the ring of differential operators 
$D_{\BC^I}=\otimes_I D_\BC$, the CDO $D^{ch}_{\BC^{I|I}}=\otimes_I D^{ch}_{\BC^{1|1}}$.

\begin{sloppypar}
As a practical matter,  what for $D^{ch}_{\BC^{1|1}}$ was a quadruple of fields, $x(z),\partial_x(z),\xi(z),\partial_{\xi}(z)$ with the OPEs recorded in 3.1, 3.2, for $D^{ch}_{\BC^{I|I}}$  is an $I$-family of quadruples, $x_i(z),\partial_{x_i}(z),\xi_i(z),\partial_{\xi_i}(z)$,
$i\in I$, which pairwise commute if the corresponding indices $i$ are distinct, and satisfy the old OPE if they are the same.
\end{sloppypar}

Whatever we know about $D^{ch}_{\BC^{1|1}}$ tends  to carry over to $D^{ch}_{\BC^{I|I}}$ immediately.  $D^{ch}_{\BC^{1|1}}$ being the
Koszul complex with differential $d^{ch}_K=\int x(z)\partial_\xi(z)\,dz$,  $D^{ch}_{\BC^{I|I}}$ is the tensor product complex with differential $d^{ch}_K= \sum_I\int x_i(z)\partial_{\xi_i}(z)\,dz$.

$D^{ch}_{\BC^{I|I}}$ is a conformal vertex algebra, with Virasoro field $L(z)=\sum_IL^{(i)}(z)$, $L^{(i)}(z)$ being the copy of the $L(z)$ from 3.3. that acts on the $i$-th tensor factor of $D^{ch}_{\BC^{I|I}}$. The $_{(1)}$-component of this field defines on $D^{ch}_{\BC^{I|I}}$
a conformal grading, and its 0-th component is the classical Koszul complex of the ring $\BC[\BC^I]$ w.r.t. the regular sequence   $ \{x_i,\, i\in I\}$, to be denoted, as usual, by $K(\BC[\BC^I], \{x_i,\, i\in I\})$, cf. 2.2.

This structure implies: the embedding of complexes

\[
K(\BC[\BC^I], \{x_i,\, i\in I\})\hookrightarrow D^{ch}_{\BC^{I|I}}
\]
is a quasiisomorphism, and 
\[
H(\BC[\BC^I], \{x_i,\, i\in I\})\stackrel{\sim}{\rightarrow} H_{d^{ch}_{K}}(D^{ch}_{\BC^{I|I}})=\BC
\]
spanned by the class of the cocycle $1=\otimes_I1$, similarly to 3.3.

\begin{sloppypar}
{\bf 4.2.}  The other half of the chiral Janus complex is again the tensor product 
$\otimes_I (D^{ch}_{\BC^{1|1}}[x_{(-1)}^{-1}]/D^{ch}_{\BC^{1|1}})$. This notation is too cumbersome to be tolerated, and for the lack of a better one we shall write
\[
(D^{ch}_{\BC^{I|I}})^d=\otimes_I (D^{ch}_{\BC^{1|1}}[x_{(-1)}^{-1}]/D^{ch}_{\BC^{1|1}})
\]
and call it the {\em vertex dual} of $D^{ch}_{\BC^{I|I}}$. The right hand side of this equality is  is indeed the dual of $D^{ch}_{\BC^{I|I}}$ in the sense of [FBZ], 10.4.6, which has already been  cited, 3.6, but since we are avoiding the discussion of this duality functor, let us take the right hand side of this equality for the definition of its left hand side.
\end{sloppypar}

The properties of $(D^{ch}_{\BC^{I|I}})^d$ are parallel to those of $D^{ch}_{\BC^{1|1}}[x_{(-1)}^{-1}]/D^{ch}_{\BC^{1|1}}$ to the same extent those of $D^{ch}_{\BC^{I|I}}$ are parallel to those of $D^{ch}_{\BC^{1|1}}$.  As a complex it is the dual of $D^{ch}_{\BC^{I|I}}$,
and so its cohomology is 1 dimensional and spanned by the class of cocycle 
\[
\prod_I\frac{\xi_{i,(-1)}}{x_{i,(-1)}},
\]
see the beginning of 3.5.

Analogously to 3.5, we define a chiral Janus, $(D^{ch}_{\BC^{I|I}})^d\stackrel{.}{\oplus} D^{ch}_{\BC^{I|I}}$, to be the direct sum
$(D^{ch}_{\BC^{I|I}})^d\oplus D^{ch}_{\BC^{I|I}}$ as a vector space (and as a $D^{ch}_{\BC^{I|I}}$-module) with differential  obtained
by tweaking $d^{ch}_K\oplus d^{ch}_K$ so as to make sure that it sends the cocycle
\[
\prod_I\frac{\xi_{i,(-1)}}{x_{i,(-1)}}\in (D^{ch}_{\BC^{I|I}})^d
\mbox{ to } 
\otimes_I1\in D^{ch}_{\BC^{I|I}},
\]
thereby making $(D^{ch}_{\BC^{I|I}})^d\stackrel{.}{\oplus} D^{ch}_{\BC^{I|I}}$ an acyclic complex.

\

{\bf 4.3. The negative (dual) chiral Koszul as a local cohomology of the positive one. }  The modules 
$(D^{ch}_{\BC^{I|I}})^d$ afford the following geometric interpretation.

The polynomial nature of $D^{ch}_{\BC^{I|I}}$ has already been   used to define localization $D^{ch}_{\BC^{1|1}}[x_{(-1)}^{-1}]$ in  3.4. One key observation made in [MSV] is that one can extend a vertex algebra structure to any localization
$D^{ch}_{\BC^{I|I}}[S^{-1}]$, $S$ being a multiplicative subset of  $\BC[\BC^I]$. This makes $D^{ch}_{\BC^{I|I}}$ into a sheaf on $\BC^I$,
to be denoted  $\underline{D}^{ch}_{\BC^{I|I}}$

{\bf 4.3.1. Lemma.}{\em  Let $U=\BC^{I}-\{0\}$, the ``principal affine space,'' $Z=\{0\}\in\BC^I$. We have

(i) For any $I$,
\[
(D^{ch}_{\BC^{I|I}})^d=H^{|I|}_Z(\BC^I, \underline{D}^{ch}_{\BC^{I|I}}).
\]

(ii) If $|I|>1$, then
\[
(D^{ch}_{\BC^{I|I}})^d=H^{|I|}_Z(\BC^I, \underline{D}^{ch}_{\BC^{I|I}})=H^{|I|-1}(U, \underline{D}^{ch}_{\BC^{I|I}}).
\]

%(ii) If $|I|=1$, then we have an exact sequence
%\[
%0\rightarrow  \underline{D}^{ch}_{\BC^{I|I}}(\BC^1)\rightarrow  \underline{D}^{ch}_{\BC^{I|I}}(\BC^1-\{0\})\rightarrow
%(D^{ch}_{\BC^{1|1}})^d\rightarrow 0,
%\]
%where $ \underline{D}^{ch}_{\BC^{I|I}}(\BC^1)\rightarrow  \underline{D}^{ch}_{\BC^{I|I}}(\BC^1-\{0\})$ is the restriction map.}
}
{\em Proof.}  Point (i) is of course well known, but we shall sketch the proof. 
Note that the point $Z=\{0\}$ is defined by the ideal $(x_i,i\in I)$. It is known, [H, Theorem 2.3], that the local cohomology
$H^{j}_Z(\BC^I, \underline{D}^{ch}_{\BC^{I|I}})$ is computed by the direct limit of the Koszul complexes:
\[
H^{j}_Z(\BC^I, \underline{D}^{ch}_{\BC^{I|I}})=\lim_{\stackrel{\rightarrow}{\{m_i\}}}H^j({D}^{ch}_{\BC^{I|I}},\{(x_{i,(-1)})^{m_i},i\in I\}).
\]
As the sequence $\{(x_{i,(-1)})^{m_i},i\in I\}$ is regular for any collection of exponents $\{m_i\}$, the cohomology is concentrated in the
$|I|$-th group. That it equals the space of ``purely singular parts,'' $\otimes_I (D^{ch}_{\BC^{1|1}}[x_{(-1)}^{-1}]/D^{ch}_{\BC^{1|1}})$, is rather easy to see.

In view of (i),  (ii) amounts to the isomorphism
\[
H^{|I|}_Z(\BC^I, \underline{D}^{ch}_{\BC^{I|I}})=H^{|I|-1}(U, \underline{D}^{ch}_{\BC^{I|I}}),
\]
which  is even better known than (i) and constitutes part of [H, Prop.2.2]. $\qed$

%It is rather clear that $D^{ch}_{\BC^{1|1}}[x_{(-1)}^{-1}]/D^{ch}_{\BC^{1|1}}\stackrel{.}{\oplus}D^{ch}_{\BC^{1|1}}$ is isomorphic to
%$D^{ch}_{\BC^{1|1}}[x_{(-1)}^{-1}]$ as a complex, but the vertex algebra structure is lost.

\

\

\centerline{\bf References}

\

[BD] A.~Beilinson, V.~Drinfeld. Chiral algebras. American
Mathematical Society Colloquium Publications, 51. American
Mathematical Society, Providence, RI, 2004. vi+375 pp. ISBN:
0-8218-3528-9

[FBZ] E.~Frenkel, D.~Ben-Zvi,  Vertex algebras and algebraic
curves. Second edition. Mathematical Surveys and Monographs, 88.
American Mathematical Society, Providence, RI, 2004;

%[CE] H.Cartan, S.Eilenberg, Homological algebra 

[G] S.I.Gelfand, Sheaves on $\BP_n$ and problems of linear algebra, Appendix to: C.Okonek, M.Schneider, 
H.Spindler, Vector bundles on complex projective spaces.

[GMS] V.Gorbounov, F.Malikov, V.Schechtman, Gerbes of chiral differential operators, II. Vertex algebroids, {\it Inv. Math.}

[H] R.Hartshorne, Local cohomology, {\it Lect. Notes in Math.} {\bf 41}.    

[KS] M.Kapranov, V.Schechtman, PROBs  and perverse sheaves II. Ran spaces and $0$-cycles with coefficients, 
arXiv:2209.02400 

%[MS] F.Malikov, V.Schechtman, Chiral de Rham complex over locally complete intersections

[MS] F.Malikov, V.Schechtman, Chiral Poincar\'e duality, {\it Math. Res. Letters}\ {\bf 6} (1999), 533-546. 

[MSV] F.Malikov, V.Schechtman, A.Vaintrob, Chiral de Rham complex, {\it Comm. Math. Phys.}

[ManS] Yu.I.Manin, V.V.Schechtman, Arrangements of hyperplanes, higher braid groups and  higher Bruhat 
orders, {\it Adv.  Studies in Pure Math.} {\bf 17} (1989), 289-308.  

[P] T.Pirashvili, On the PROP corresponding to bialgebras, arXiv:math/0110014;  {\it Cah. Topol G\'eom. Diff\'er. Cat\'eg.} {\bf 43} (2002), 221-239.

[S] V.Schechtman, PROBs and sheaves, talk on a joint work with M.Kapranov, Weizmann Institute of Science, 
December 2021, 
\newline https://www.math.univ-toulouse.fr/~schechtman/probs-weizmann.pdf

\

\

F.M.: Dept. of Mathematics, University of Southern California, Los Angeles, CA 90089, USA;
fmalikov@math.usc.edu

\

V.S.: Institut de Math\'ematiques de Toulouse, Universit\'e Paul Sabatier, 118 Route de Narbonne, 31062 Toulouse, 
France; Kavli IPMU, University of Tokyo, 5-1-5 Kashiwa-no-ha, Kashiwa-shi, Chiba, 277-8583 Japan; schechtman@math.univ-toulouse.fr

\

\

\end{document}